\documentclass[11pt]{article}
\usepackage{hyperref}
\pdfoutput=1
\usepackage[english]{babel}
\usepackage{amsmath,amsthm}
\usepackage{subeqnarray}
\usepackage{amsfonts}
\usepackage{graphicx}
\usepackage[onehalfspacing]{setspace}
\usepackage[left=1.45in, top=1.25in, right=0.75in, bottom=1.15in]{geometry}
\newtheorem{thm}{Theorem}[section]

\theoremstyle{definition}

\theoremstyle{remark}
\newtheorem{remark}[thm]{Remark}
\numberwithin{equation}{section}

\newcommand{\set}[2]{\left\{{#1}\,:~{#2}\right\}}
\newcommand {\average}[1] {\mbox{$\left\{\!\!\left\{ #1 \right\}\!\!\right\}$}}
\newcommand {\jump}[1] {\mbox{$\left[\!\left[ #1 \right]\!\right]$}}

\begin{document}

\title{Space-Time Discontinuous Galerkin Solution of Convection Dominated Optimal Control Problems}
\author{Tu\u{g}ba Akman \thanks{Electronic address: \texttt{takman@metu.edu.tr}} and B\"ulent Karas\"ozen \thanks{Electronic address: \texttt{bulent@metu.edu.tr}}\\
\\\vspace{6pt} Department of Mathematics and Institute of Applied Mathematics \\
Middle East Technical University, 06800 Ankara, Turkey}

\maketitle

\begin{abstract}
  In this paper, a space-time discontinuous Galerkin finite element method for distributed optimal control problems governed by unsteady diffusion-convection-reaction equations with control constraints  is studied. Time discretization is performed by discontinuous Galerkin method with piecewise constant and linear polynomials, while symmetric interior penalty Galerkin with upwinding is used for space discretization. The numerical results presented confirm the theoretically observed convergence rates.
\end{abstract}


\section{Introduction} \label{introduction}
Optimal control problems (OCPs) governed by diffusion-convection-reaction equations arise in environmental control problems, optimal control of fluid flow and in many other applications. It is well known that the standard Galerkin finite element discretization causes nonphysical oscillating solutions when convection dominates. Stable and accurate numerical solutions can be achieved by various effective stabilization techniques such as the streamline upwind/Petrov Galerkin (SUPG) finite element method \cite{SSCollis_MHeinkenschloss_2002a}, the local projection stabilization \cite{RBecker_BVexler_2007a}, the edge stabilization \cite{MHinze_NYan_ZZhou_2009a}. Recently, discontinuous Galerkin (dG) methods gain importance due to their better convergence behaviour, local mass conservation, flexibility in approximating rough solutions on complicated meshes, mesh adaptation and weak imposition of the boundary conditions in OPCs, see, e.g., \cite{DLeykekhman_2012b,DLeykekhman_MHeinkenschloss_2012a,HYucel_MHeinkenschloss_BKarasozen_2012b,HYucel_BKarasozen_2013}.\\

In the recent years much effort has been spent  on parabolic OCPs (see for example \cite{Apel_Flaig_2012_a, DMeidner_BVexler_2008b}). There are few publications dealing with OCPs governed by nonstationary diffusion-convection-reaction equation. The local DG approximation of the OCP which is discretized by backward Euler in time is studied in \cite{ZZhou_NYan_2010a} and a priori error estimates for semi-discrete OCP is provided in \cite{TSun_2010a}. In \cite{HFu_2010a,HFu_HRui_2009a}, the characteristic finite element solution of the OCP is discussed and numerical results are provided. A priori error estimates for discontinuous Galerkin time discretization for unconstrained parabolic OCPs is proposed in \cite{CKonstantinos_2007a}. Crank-Nicolson time discretization is applied to OCP of diffusion-convection equation in \cite{EBurman_2011a}. To the best of our knowledge, this is the first study on space-time dG discretization of unsteady OCPs governed by convection-diffusion-reaction equations. \\

There are two different approaches for solving OCPs: \textit{optimize-then-discretize (OD)} and \textit{discretize-then-optimize (DO)}. In the \emph{OD} approach, first the infinite dimensional optimality system is derived containing state and adjoint equation and the variational inequality. Then, the optimality system is discretized by using a suitable discretization method in space and time. In \emph{DO} approach, the infinite dimensional OCP is discretized and then the finite-dimensional optimality system is derived. The
\emph{DO}  and  \emph{DO}   approaches do not commute in general for OCPs governed by diffusion-convection-reaction equation \cite{SSCollis_MHeinkenschloss_2002a}. However, commutativity is achieved in the case of SIPG discretization for steady state problems \cite{DLeykekhman_2012b, HYucel_MHeinkenschloss_BKarasozen_2012b}. For dG time discretization, we show that \emph{OD} and \emph{DO} approaches commute also for time-dependent problems. \\

In this paper, we solve the OCP governed by diffusion-convection-reaction equation with control constraints by applying symmetric interior penalty Galerkin (SIPG) method in space and discontinuous Galerkin (dG) discretization in time \cite{KChrysafinos_NJWalkington_2006b, KEriksson_CJohnson_VThomee_1985, MFeistauer_VKucera_KNajzar_JProkopovza_2011a, VThomee_1997, MVlasak_VDolejsi_JHajek_2010a}.   In the study of Konstantinos \cite{CKonstantinos_2007a}, a priori error estimates for continuous in space and discontinuous in time Galerkin discretization for unconstrained parabolic OCPs are derived by decoupling the optimality system. In \cite{HFu_HRui_2009a}, error analysis concerning the characteristic finite element solution of the OCP with control constraints is discussed.
Optimal order of convergence rates for the space-time discretization is confirmed on two numerical examples. Additionally we give numerical results for Crank-Nicolson method and compare them with the DG in time discretization.\\

The rest of the paper is organized as follows. In Section \ref{ocp}, we define the model problem and then derive the optimality system. In Section \ref{dg}, discontinuous Galerkin discretization and the semi-discrete optimality system follow. In Section \ref{time}, space-time dG methods and state the fully discrete optimality system are presented.
In Section \ref{numerical}, numerical results are shown in order to discover the performance of the suggested method. The paper ends with some conclusions.

\section{The Optimal Control Problem}\label{ocp}
We consider the following distributed optimal control problem governed by the unsteady diffusion-convection-reaction equation with control constraints
\begin{subeqnarray}\label{s1}
\underset{u \in U_{ad}}   {\hbox{ minimize }}\; J(y,u):= \frac{1}{2}\int_{0}^{T} &\big(&\left\|y-y_{d}\right\|^{2}_{L^{2}(\Omega)} \ +
                    \alpha \left\|u\right\|^{2}_{L^{2}(\Omega)} \big) \; dt,\\
\hbox{ subject to } \partial_t y-\epsilon \Delta y+\beta \cdot\nabla y+r y &=&  f+u \quad (x, t) \in \Omega \times (0,T],   \\
  y(x,t)&=& 0  \quad \quad \; \; (x, t) \in  \partial \Omega \times [0,T],   \\
 y(x,0) &=& y_{0}(x) \qquad \; \; x \in \Omega,
\end{subeqnarray}
where the admissible space of control constraints  is given by
\begin{align}\label{s2_1}
U_{ad} = \{ u \in L^{2}(0,T;U): u_a \le u \le u_b, \hbox{ a.e. in } \Omega \times (0,T] \}
\end{align}
with the constant bounds $u_a, u_b \in L^{\infty}(\Omega)$, i.e., $u_a <u_b$. We take $\Omega$ as a bounded open convex domain in $\mathbb{R}^{2}$ with Lipschitz boundary $\partial \Omega$ and $I = (0, T]$ as time interval. The source function and the desired state are denoted by $f \in L^{2}(0,T;L^{2}(\Omega))$ and $y_d \in L^{2}(0,T;L^{2}(\Omega))$, respectively. The initial condition is also defined as $y_0(x) \in H_{0}^1(\Omega)$. The diffusion and reaction coefficients are  $\epsilon >0$ and $r \in L^{\infty}(\Omega)$, respectively. The velocity field $\beta \in (W^{1,\infty}(\Omega))^2$  satisfies the incompressibility condition, i.e. $\nabla \cdot \beta =0$. Furthermore, we assume the existence of constant $c_{0} \equiv c_{0}(x)$ such that $r(x) - \frac{1}{2} \nabla \cdot \beta(x) \geq c_{0}$ a.e. $\in \Omega$ so that the well-posedness of the optimal control problem (\ref{s1}) is guaranteed. The  trial and test spaces are
\[
Y=V=H^1_0(\Omega), \quad  \forall t \in (0,T].
\]
For $ (y,u) \in  Y \times U_{ad}$, the variational formulation corresponding to (\ref{s1}) can be written
\begin{subeqnarray}\label{s3}
 \underset {u \in U_{ad}}  {\hbox{ minimize }} \;
  J(y,u):= \frac{1}{2} \int_{0}^{T} &\big(& \|y-y_{d}\|^{2}_{L^2(\Omega)} + \alpha \| u\|^{2}_{L^2(\Omega)}\big)\; dt     \\
  \hbox{subject to} \; (\partial_t y,v)+ a(y,v) &=& (f + u,v), \; \forall v \in V, \; t \in I, \\
                                         y(x,0) &=& y_0, \qquad \quad \;\; x \in \Omega,
\end{subeqnarray}
with $a(y,x) = \int_{\Omega} (\epsilon \nabla y \cdot \nabla v + \beta \cdot \nabla y v + r y v )dx$ and $(u, v) =  \int_{\Omega} u v dx$. Differentiating  the Lagrangian
\begin{eqnarray*}
\mathcal{L}(y, u, p)
&=& \frac{1}{2} \int_{0}^{T} \big( \|y-y_{d}\|^{2}_{L^2(\Omega)} + \alpha \| u\|^{2}_{L^2(\Omega)}\big)\; dt  \\
&+& \left\lbrace  (\partial_t y,p)+ a(y,p) - (f + u,p) \right\rbrace + (y(x,0) - y_{0}, p(x,0)).
\end{eqnarray*}
 with respect to $p, y, u$, we obtain the optimality system
\begin{subeqnarray}\label{s4}
(\partial_t y,v)+a(y,v) & = & (f + u,v), \qquad  \; \; y(x,0)=y_0, \label{s4a} \\
 -(\partial_t p,\psi)+a(\psi,p) & = & -(y-y_{d},\psi), \quad p(x,T)=0, \label{s4b} \\
\int_{0}^{T} (\alpha u-p, w -u) dt & \geq 0 & , \quad \forall w \in U_{ad}. \label{s4c}
\end{subeqnarray}
It is well known that the functions $(y,u) \in H^{1}(0,T;L^{2}(\Omega)) \cap L^{2}(0,T;Y) \times U_{ad}$ solve (\ref{s3}) if and only if there is an adjoint $p \in H^{1}(0,T;L^{2}(\Omega)) \cap L^{2}(0,T;Y)$ such that $(y,u,p)$  is the unique solution of the optimality system (\ref{s4}) \cite{JLLions_1971,FTroltzsch_2010a}.
\section{Discontinuous Galerkin Discretization}\label{dg}
Let $\{ \mathcal{T}_h\}_h$ be a family of shape regular meshes such that $\overline{\Omega} = \cup_{K \in \mathcal{T}_h} \overline{K}$, $K_i \cap K_j = \emptyset$ for $K_i, K_j \in \mathcal{T}_h$, $i \not= j$. The diameter of an element $K$ and the length of an edge $E$ are denoted by $h_{K}$ and $h_E$, respectively. In addition, the maximum value of element diameter is denoted by $h=\max \limits_{K \in \mathcal{T}_h} h_{K}$. We only consider  discontinuous piecewise  linear finite element spaces to
define the discrete state and control spaces
\begin{eqnarray}\label{Uspace}
V_h = Y_h = U_{h}^{ad} &=& \set{y \in L^2(\Omega)}{ y\mid_{K}\in \mathbb{P}^1(K) \quad \forall K \in \mathcal{T}_h}\\
\widetilde{U}_h &\subset& (Y_{h} \cap U_{ad}),
\end{eqnarray}
\noindent respectively. Here, $\mathbb{P}^1(K)$ denotes the set of all polynomials on $ K \in \mathcal{T}_h$ of degree at most $1$.

\begin{remark}
When the state equation (\ref{s1}) contains nonhomogeneous Dirichlet boundary conditions, $Y_{h}$  and  $V_{h}$ can still be taken as the same space due to the weak treatment of boundary conditions in dG methods (see for example \cite{DLeykekhman_MHeinkenschloss_2012a}).
\end{remark}

We split the  set of all edges $\mathcal{E}_h$ into the set $\mathcal{E}^0_h$ of interior edges  and the set $\mathcal{E}^{\partial}_h$ of boundary edges so that $\mathcal{E}_h=\mathcal{E}^{\partial}_h\cup \mathcal{E}^{0}_h$. Let $\mathbf{n}$ denote the unit outward normal to $\partial \Omega$. We define the inflow boundary
\[
          \Gamma^- = \set{x \in \partial \Omega}{ \beta\cdot \mathbf{n}(x) < 0}
\]
and the outflow boundary $\Gamma^+ = \partial \Omega \setminus  \Gamma^-$. The boundary edges are decomposed into edges
$\mathcal{E}^{-}_h = \set{ E \in \mathcal{E}^{\partial}_h}{ E \subset  \Gamma^- }$ that correspond to inflow boundary and edges
$\mathcal{E}^{+}_h = \mathcal{E}^\partial_h \setminus \mathcal{E}^{-}_h$ that  correspond to outflow boundary. The inflow and outflow boundaries of an element $K \in \mathcal{T}_h$ are defined by
\begin{equation*}
\partial K^- =\set{x \in \partial K}{\beta \cdot \mathbf{n}_{K}(x) <0}, \quad \partial K^{+} = \partial K \setminus \partial K^{-},
\end{equation*}
where $\mathbf{n}_{K}$ is the unit normal vector on the boundary $\partial K$ of an element $K$.

Let the edge $E$ be a common edge for two elements $K$ and $K^e$. For a piecewise continuous scalar function $y$,
there are two traces of $y$ along $E$, denoted by $y|_E$ from interior of  $K$ and $y^e|_E$ from interior of  $K^e$. Then, the jump and average of $y$ across the edge $E$ are defined by:
\begin{equation}
\jump{y} =y|_E\mathbf{n}_{K}+y^e|_E\mathbf{n}_{K^e}, \quad
\average{y}=\frac{1}{2}\big( y|_E+y^e|_E \big).
\end{equation}

Similarly, for a piecewise continuous vector field $\nabla y$, the jump and average across an edge $E$ are given by
\begin{equation}
\jump{\nabla y} =\nabla y|_E \cdot \mathbf{n}_{K}+\nabla y^e|_E \cdot \mathbf{n}_{K^e}, \quad
\average{\nabla y}=\frac{1}{2}\big(\nabla y|_E+\nabla y^e|_E \big).
\end{equation}

For a boundary edge $E \in K \cap \Gamma$, we set $\average{\nabla y}=\nabla y$ and $\jump{y}=y\mathbf{n}$
where $\mathbf{n}$ is the outward normal unit vector on $\Gamma$.

We can now give dG discretizations of the state equation (\ref{s1}) in space  for fixed control $u$.  The dG method proposed here is based on the upwind discretization of the convection term and on the SIPG discretization  of the diffusion term \cite{DSchotzau_LZhu_2009a}. This leads to the following (bi-)linear forms applied to $y_{h} \in H^1(0,T;Y_h)$ for $\forall t \in (0,T]$
\begin{equation}
(\partial_t y_h, v_h) + a_h^{s}(y_h,v_h)+b_h(u_h,v_h)=(f,v_h)  \quad \forall v_h \in V_h, \quad  t \in (0,T],
\end{equation}
where
\begin{eqnarray}\label{s5}
a^{d}(y,v) &=& \sum \limits_{K \in \mathcal{T}_h} \int \limits_{K} \epsilon \nabla y \cdot  \nabla v \; dx \nonumber \\
           &-& \sum \limits_{ E \in \mathcal{E}_h} \int \limits_E \big( \average{\epsilon \nabla y} \cdot \jump{v}
           + \average{\epsilon \nabla v} \cdot \jump{y}
           - \overbrace{\frac{\sigma \epsilon}{h_E}  \jump{y} \cdot \jump{v} \big) \; ds}^{J_{\sigma}(y, v)}
\end{eqnarray}
and
\begin{eqnarray}\label{s5}
a_h^{s}(y,v) &=& a^{d}(y,v)
         + \sum \limits_{K \in \mathcal{T}_h} \int \limits_{K} \big( \beta \cdot \nabla y v + r y v \big) \; dx \nonumber\\
         &+& \sum \limits_{K \in \mathcal{T}_h}\; \int \limits_{\partial K^{-} \backslash \Gamma^{-}} \beta \cdot \mathbf{n} (y^e-y)v \; ds
         - \sum \limits_{K \in \mathcal{T}_h} \; \int \limits_{\partial K^{-} \cap \Gamma^{-}} \beta \cdot \mathbf{n} y v  \; ds, \\
b_h(u, v) &=& - \sum \limits_{K \in \mathcal{T}_h} \int \limits_{K} uv \; dx
\end{eqnarray}
with a constant interior penalty parameter $\sigma >0$. We choose $\sigma$ to be sufficiently large, independent of the mesh size $h$ and the diffusion coefficient $\epsilon$ to ensure the stability of the dG discretization as described in \cite[Sec.~2.7.1]{BRiviere_2008a} with a lower bound depending  only on the polynomial degree. Large penalty parameters  decrease the jumps across element interfaces, which can affect the numerical approximation \cite{DNArnold_FBrezzi_BCockburn_LDMarini_2001a}.

\subsection{Semi-discrete Formulation of The Optimal Control Problem}\label{semi}
\noindent Let $f_h, y_h^d$ and $y_h^0$ be approximations of the source function $f$, the desired state function $y_d$ and initial condition $y_0$, respectively. Then, the semi-discrete approximation of the optimal control problem (\ref{s4}) can be defined as follows:
\begin{subeqnarray}\label{s6}
 \underset{u_h \in U_h^{ad}} {\hbox{ minimize }} \int_{0}^{T} \big(\frac{1}{2} \sum \limits_{K \in \mathcal{T}_h} \|y_h-y_{h}^d\|^{2}_{L^{2}(K)} &+& \frac{\alpha}{2} \sum \limits_{K \in \mathcal{T}_h} \| u_h\|^{2}_{L^{2}(K)} \big) \; dt,   \\
  \hbox{subject to } (\partial_t y_h, v_h) + a_h^{s}(y_h,v_h)+b_h(u_h,v_h) &=& (f_h,v_h), \label{s6_b}\\
                        y_h(x,0)&=&y_h^0, \nonumber \\
                        t \in (0,T], v_h \in V_h, (y_h,u_h) &\in& Y_h \times U_h^{ad}. \nonumber
\end{subeqnarray}

The semi-discrete optimality system is written as follows:
\begin{eqnarray}\label{s7}
(\partial_t y_h, v_h) + a_{h}^{s}(y_{h},v_{h})+b(u_{h},v_{h})&=&(f_h,v_{h}), \qquad \qquad \; y_h(x,0)=y_h^0, \label{s7a} \nonumber \\
-(\partial_t p_h, \psi_{h}) + a_{h}^{a}(p_{h}, \psi_{h})&=&-(y_{h}-y_h^{d},\psi_{h}), \quad p_h(x,T)=0,\label{s7b}  \\
\int_{0}^{T} (\alpha u_{h}-p_{h}, w_{h} -u_{h})\; dt  &\geq& 0, \quad \forall w_{h} \in U^{ad}_{h}, \nonumber \label{s7c}
\end{eqnarray}

where
\begin{eqnarray*}\label{s55}
a_h^{a}(p,\psi) &=& \sum \limits_{K \in \mathcal{T}_h} \int \limits_{K} \epsilon \nabla p \cdot  \nabla \psi \; dx \nonumber \\
             &-& \sum \limits_{ E \in \mathcal{E}_h} \int \limits_E \big( \average{\epsilon \nabla p} \cdot \jump{\psi}
             + \average{\epsilon \nabla \psi} \cdot \jump{p}
             - \frac{\sigma \epsilon}{h_E}  \jump{p} \cdot \jump{\psi} \big) \; ds\nonumber \\
             &+& \sum \limits_{K \in \mathcal{T}_h} \int \limits_{K} \big( -\beta \cdot \nabla p \psi + r p \psi \big) \; dx \nonumber\\
             &-& \sum \limits_{K \in \mathcal{T}_h}\; \int \limits_{\partial K^{+} \backslash \Gamma^{+}} \beta \cdot \mathbf{n} (p^e-p)\psi \; ds
              + \sum \limits_{K \in \mathcal{T}_h} \; \int \limits_{\partial K^{+} \cap \Gamma^{+}} \beta \cdot \mathbf{n} p \psi  \; ds.
\end{eqnarray*}

\section{Time Discretization of The Optimal Control Problem}\label{time}
In this section, we derive the fully-discrete optimality system, by using $\theta$-method and discontinuous Galerkin method. We compare the resulting optimality systems based on two approaches, i.e. \emph{optimize-then-discretize (OD)} and \emph{discretize-then-optimize (DO)}. Let $0 = t_0 < t_1 < \cdot < t_{N_{T}} = T $ be a subdivison of $I=(0,T)$ with time intervals $I_m = (t_{m-1}, t_{m}]$ and time steps $k_{m} = t_{m} - t_{m-1}$ for $m=1, \ldots, N_{T}$ and $k = \max_{1 \leq m \leq N_{T}} k_{m}$.

\subsection{Time Discretization Using $\theta$-method}\label{theta}
We start with \emph{OD} approach by discretizing the semi-discrete optimality system (\ref{s7}) using $\theta$-method.
\begin{eqnarray}\label{OD_ocp_theta}
& & (y_{h,m+1}-y_{h,m},v) + k a_{h}^{s}((1-\theta) y_{h,m}+\theta y_{h,m+1},v) =  \nonumber\\
& & k((1-\theta) f_{h,m} + \theta f_{h,m+1}) + k((1-\theta) u_{h,m}+\theta u_{h,m+1},v), \; m=0, \cdots, N-1,  \nonumber \\
& &  y_{h,0}(x,0)=y_{0} \nonumber \\
& & (p_{h,m} - p_{h,m+1},q) + k a_{h}^{a}(\theta p_{h,m}+ (1-\theta)p_{h,m+1},q) =   \\
&&  -k\left(\theta(y_{h,m} - y^{d}_{h,m},q) + (1-\theta)(y_{h,m+1}-y^{d}_{h,m+1},q)\right),\;  m = N-1, \cdots, 0,\nonumber\\
& & p_{h,N} = 0, \nonumber \\
& & (\alpha u_{h,m}-p_{h,m},w-u_{h,m}) \geq 0 \quad m=0,1, \ldots, N. \nonumber
\end{eqnarray}

We proceed with \emph{DO} approach. To do this, we approximate the first part of the cost functional by the rectangle rule, the second part of it by the trapezoidal rule and discretize the state equation using $\theta$-method in time. We use the rectangle rule to approximate the first part so that the value of the adjoint at the final time becomes zero as in \cite{MStoll_AWathen_2010a}.
\begin{eqnarray*}
& & \underset{u_\delta \in \widetilde{U}_{h}^{k}}  {\hbox{ minimize }}
 \frac{k}{2}\sum \limits_{m=0}^{N-1}(y_{h,m}-y^{d}_{h,m})^{T} M(y_{h,m}-y^{d}_{h,m}) \\
& &+ \alpha \frac{k}{2} \left(\frac{1}{2}u_{h,0}^{T}M u_{h,0}
+ \sum \limits_{m=1}^{N-1}u_{h,m}^{T}M u_{h,m}
+ \frac{1}{2} u_{h,N}^{T}M u_{h,N}\right) \nonumber \\
& & \hbox{ subject to } \\
& & (y_{h,m+1} - y_{h,m},v) + ka_{h}^{s}((1-\theta) y_{h,m}+\theta y_{h,m+1},v) =\\
& &  k((1-\theta) f_{h,m} + \theta f_{h,m+1})+k((1-\theta) u_{h,m}+\theta u_{h,m+1},v),
m=0, \cdots, N-1, \\
& &  (y_{h,0},v) = (y_{0},v),
\end{eqnarray*}
where $M$ is the mass matrix.\\

Now we construct the discrete Lagrangian
\begin{eqnarray}\label{Lagrangian_theta}
&\mathcal{L}&(y_{h,1}, \ldots, y_{h,N}, p_{h,0}, \ldots, p_{h,N}, u_{h,0}, \ldots, u_{h,N}) \nonumber\\
&=& \frac{k}{2}\sum \limits_{m=0}^{N-1}(y_{h,m}-y^{d}_{h,m})^{T} M(y_{h,m}-y^{d}_{h,m})\nonumber\\
&+& \alpha \frac{k}{2} \left(\frac{1}{2}u_{h,0}^{T}M u_{h,0}
+ \sum \limits_{m=1}^{N-1}u_{h,m}^{T}M u_{h,m}
+ \frac{1}{2} u_{h,N}^{T}M u_{h,N}\right) + (y_{h,0}-y_{0},p_{h,0})\nonumber\\
&+& \sum \limits_{m=0}^{N-1}((y_{h,m+1}-y_{h,m},p_{h,m+1}) + ka_{h}^{s}((1-\theta) y_{h,m}+\theta y_{h,m+1},p_{h,m+1}) \nonumber\\
&-& k((1-\theta) f_{h,m} + \theta f_{h,m+1})+k((1-\theta) u_{h,m}+\theta u_{h,m+1},p_{h,m+1})).
\end{eqnarray}

By differentiating Lagrangian (\ref{Lagrangian_theta}), we derive the fully-discrete optimality system
\begin{eqnarray}\label{DO_ocp_theta}
& & (y_{h,m+1} - y_{h,m},v) + k a_{h}^{s}((1-\theta) y_{h,m}+\theta y_{h,m+1},v) = \nonumber\\
& & k((1-\theta) f_{h,m} + \theta f_{h,m+1}) + k((1-\theta) u_{h,m}+\theta u_{h,m+1},v),\; m = 0, \cdots, N-1 \nonumber \label{DO_ocp_theta_a1}  \\
& & y_{h,0}(x,0)=y_{0} \nonumber \label{DO_ocp_theta_a2}\\
& & (q,p_{h,N}) + ka^{s}_{h}(q,\theta p_{h,N}) = 0,  \nonumber \label{DO_discerete_ocp_theta_b1}\\
& & (p_{h,m} - p_{h,m+1},q) + k a_{h}^{s}(q,\theta p_{h,m}+(1-\theta)p_{h,m+1}) =\\
& & -k(y_{h,m} - y^{d}_{h,m},q),\;  m = N-1,\ldots,1, \nonumber \label{DO_discerete_ocp_theta_b2}\\
& & (q,p_{h,0}-p_{h,1})+ ka^{s}_{h}(q,(1-\theta)p_{h,1}) = -k(y_{h,0} - y^{d}_{h,0},q),  \nonumber \label{DO_discerete_ocp_theta_b3}\\
& & (\frac{\alpha}{2}u_{h,0} - (1-\theta)p_{h,1},w-u_{h,0}) \geq 0,  \nonumber  \label{DO_discerete_ocp_theta_c1} \\
& & (\alpha u_{h,m} - (\theta p_{h,m}+(1-\theta)p_{h,m+1}),w-u_{h,m}) \geq 0, \quad m=1, \ldots, N-1, \nonumber \label{DO_discerete_ocp_theta_c2} \\
& & (\frac{\alpha}{2}u_{h,N} - \theta p_{h,N},w-u_{h,N}) \geq 0. \nonumber \label{DO_discerete_ocp_theta_c3}
\end{eqnarray}

In the case of backward Euler method ($\theta = 1$), the value $u_{h,0}$ is not needed as we observe from (\ref{DO_ocp_theta_a1}). As we mentioned before, the approximation of the first integral in the cost functional by using the rectangle rule leads to $p_{h,N} = 0$, $u_{h,N} = 0$, as we see from (\ref{DO_discerete_ocp_theta_c2}). For the SIPG we obtain $a_{h}^{s}(\psi_{\delta},p_{\delta}) = a_{h}^{a}(p_{\delta},\psi_{\delta})$ \cite{HYucel_MHeinkenschloss_BKarasozen_2012b}
and therefore  (\ref{OD_ocp_theta}) and (\ref{DO_ocp_theta}) gives the same variational formulation.

In the case of Crank-Nicolson method ($\theta = 1/2$), we observe that some differences occur in the adjoint equation. In (\ref{OD_ocp_theta}), the right-hand side of the adjoint equation is evaluated at two successive points, while it is evaluated at just one point in (\ref{DO_ocp_theta}). Additional differences are seen in the variational inequalities (\ref{OD_ocp_theta}) and (\ref{DO_ocp_theta}), too. Thus, \emph{OD} and \emph{DO} approaches lead to different weak formulations. In~\cite{Apel_Flaig_2012_a}, the optimal control of the heat equation is concerned by applying continuous Galerkin discretization. For \emph{DO} approach, the cost functional is discretized by using the midpoint rule. On the other hand, for \emph{OD} approach, the semi-discrete state equation is discretized by using the midpoint rule and a variation of the trapezoidal rule is applied to the semi-discrete adjoint equation to obtain the fully discrete optimality system. Then \emph{OD} and \emph{DO} approaches commute.

\subsection{Time Discretization Using Discontinuous Galerkin Method}\label{dgintime}
We derive the fully discrete optimality system by employing discontinuous Galerkin time discretization to the semi-discrete optimality system (\ref{s7}). We define the space-time finite element space of piecewise discontinuous functions for state and control as
\begin{eqnarray*}
V_{h}^{k,q} = Y_{h}^{k,q} &=& \set{ v \in L^{2}(0,T; L^{2}(\Omega))}{v |_{I_{m}} = \sum_{s=0}^{q} t^{s} \phi_{s}, t \in I_{m}, \phi_{s} \in V_{h}, m = 1, \ldots, N},\\
\widetilde{U}_h^{k, q} &\subset& (Y_{h}^{k,q} \cap U_{ad}).
\end{eqnarray*}

We define the temporal jump of $v \in V_{h}^{k,q}$ as $[v]_{m} = v_{+}^{m} - v_{-}^{m}$, where $w_{\pm}^{m} = \lim\limits_{\varepsilon \rightarrow 0 \pm} v(t_{m} + \varepsilon)$. Let $f_{\delta}$ and $y_{\delta}^d$ be approximations of the source function $f$ and the desired state function $y^d$ on each interval $I_m$. Then, the fully-discrete optimal control problem is written as
\begin{subeqnarray} \label{fully_discerete_cost}
\underset {u_{\delta} \in \widetilde{U}_h^{k, q}}   {\hbox{ minimize }}
\frac{1}{2} \int_{0}^{T} \sum \limits_{K \in \mathcal{T}_h}&& \big( \|y_{\delta}-y_{\delta}^{d}\|^{2}_{L^{2}(K)} + \alpha \| u_{\delta}\|^{2}_{L^{2}(K)} \big)  dt, \\
\hbox{ subject to }
\int_{0}^{T}\big((\partial_t y_{\delta}, v_{\delta}) &+& a_{h}^{s}(y_{\delta},v_{\delta})\big) dt + \sum \limits_{m=1}^{N_T}([y_{\delta}]_{m-1}, v^{m-1}_{\delta, +} ) \nonumber \\
&=&  \int_{0}^{T} (f_\delta + u_{\delta}, v_{\delta})dt, \qquad y_{\delta, 0}^{-} = (y_0)_{\delta}.
\end{subeqnarray}

The OCP (\ref{fully_discerete_cost}) has a unique solution $(y_{\delta}, u_{\delta})$ and that pair $(y_{\delta}, u_{\delta}) \in V_{h}^{k, q} \times \widetilde{U}_h^{k, q}$ is the solution of (\ref{fully_discerete_cost}) if and only if there is an adjoint $p_{\delta} \in V_{h}^{k, q}$ such that $(y_{\delta}, u_{\delta}, p_{\delta}) \in  V_{h}^{k, q} \times \widetilde{U}_h^{k, q} \times V_{h}^{k, q}$  is the unique solution of the fully-discrete optimality system
\begin{subeqnarray}\label{fully_discerete_ocp}
\int_{0}^{T}\big((\partial_t y_{\delta}, v_{\delta}) + a_{h}^{s}(y_{\delta},v_{\delta})\big) dt + \sum \limits_{m=1}^{N_T}([y_{\delta}]_{m-1}, v^{m-1}_{\delta, +} )
&=& \int_{0}^{T}(f_\delta + u_{\delta}, v_{\delta})\; dt, \notag\\
y_{\delta, 0}^{-} &=& (y_0)_{\delta}, \label{fully_discerete_ocp_a} \\
\int_{0}^{T}\big(-(\partial_t p_{\delta}, \psi_{\delta}) + a_{h}^{a}(p_{\delta},\psi_{\delta})\big)  dt -   \sum \limits_{m=1}^{N_T}([p_{\delta}]_{m}, \psi^{m}_{\delta, -} )
&=& -  \int_{0}^{T}(y_{\delta}-y_{\delta}^{d},\psi_{\delta})\; dt, \notag \\
p_{\delta, N}^{+} &=& 0,\label{fully_discerete_ocp_b}  \\
\int_{0}^{T} (\alpha u_{\delta} - p_{\delta}, w_{\delta} - u_{\delta})\; dt  &\geq& 0 \quad  \forall w_{\delta} \in \widetilde{U}_h^{k, q}. \label{fully_discerete_ocp_c}
\end{subeqnarray}
We note that (\ref{fully_discerete_ocp}) is obtained by discretizing (\ref{s4}), that is, we employ OD approach.

Finally, we define the auxiliary problem which is needed for a priori error analysis
\begin{equation}
(J'_{\delta}(u), v - u) = \int_{0}^{T}(\alpha u - p_{\delta}^{u}, v - u) dt \label{aux_cost_func},
\end{equation}
subject to
\begin{subeqnarray}\label{aux_fully_discerete_ocp}
\int_{0}^{T}\left(  (\partial_t y_{\delta}^{u}, v_{\delta}) + a_{h}^{s}(y_{\delta}^{u},v_{\delta})\right) dt +  \sum \limits_{m=1}^{N_T}([y_{\delta}^{u}]_{m-1}, v^{m-1}_{\delta, +} )
&=& \int_{0}^{T}(f_\delta + u, v_{\delta})\; dt, \notag\\
y_{\delta, 0}^{u -} &=& (y_0)_{\delta}, \label{aux_fully_discerete_ocp_a} \\
 \int_{0}^{T}\left(-(\partial_t p_{\delta}^{u}, \psi_{\delta}) + a_{h}^{a}(p_{\delta}^{u},\psi_{\delta})\right)  dt -   \sum \limits_{m=1}^{N_T}([p_{\delta}^{u}]_{m}, \psi^{m}_{\delta, -} )
&=& -  \int_{0}^{T}(y_{\delta}^{u}-y_{\delta}^{d},\psi_{\delta})\; dt, \notag \\
p_{\delta, N}^{u +} &=& 0.\label{aux_fully_discerete_ocp_b}
\end{subeqnarray}

\subsection{Commutativity Properties of Space-Time dG Method}\label{commutativity}
In the case of time-dependent OCP, the difference between the optimality system arising from \emph{OD} and \emph{DO} is caused by nonsymmetric nature of the bilinear form or the inconsistency of the final condition of the adjoint equation with the optimality system. In the DO approach, we construct the discrete Lagrangian
\begin{eqnarray*}
\mathcal{L}(y_{\delta}, u_{\delta}, p_{\delta})
&=& \frac{1}{2} \int_{0}^{T} \left(  \sum \limits_{K \in \mathcal{T}_h} \big( \|y_{\delta}-y_{\delta}^{d}\|^{2}_{L^{2}(K)} + \alpha \| u_{\delta}\|^{2}_{L^{2}(K)} \big)\right)  dt\\
&+& \sum_{m=1}^{N_T} \big(\int_{I_m}\left(  (\partial_t y_{\delta}, p_{\delta}) + a_{h}^{s}(y_{\delta},p_{\delta})\right) dt +  ([y_{\delta}]_{m-1}, p^{m-1}_{\delta, +})\big) \\
&-& \sum_{m=1}^{N_T}  \int_{I_m} (f_\delta + u_{\delta}, p_{\delta}) dt\big)
+ ((y_0)_{\delta} - y_{\delta, 0}^{-}, p_{\delta, 0}^{-}).
\end{eqnarray*}
Differentiating $\mathcal{L}$ with respect to $y_{\delta}$ and applying integration by parts, we obtain
\begin{eqnarray}\label{adjoint_DO1}
\sum_{m=1}^{N_T} \int \limits_{I_m} \big( \psi_{\delta}, -\partial_t p_{\delta})
&+& a_{h}^{s}(\psi_{\delta}, p_{\delta})\big)  dt
+ \sum \limits_{m=1}^{N_T-1}(\psi^{-}_{\delta, m}, -[p_{\delta}]_{m} ) + (q_{\delta, N_{T}}^{-}, p_{\delta, N_T}^{-})\notag\\
&=& - \sum_{m=1}^{N_T} \int_{I_m}(y_{\delta}-y_{\delta}^{d},\psi_{\delta})\; dt, \quad \forall \psi_{\delta} \in V_{h}^{k, q}.
\end{eqnarray}
Now, we add and subtract $(\psi_{\delta, N_{T}}^{-}, p_{\delta, N_T}^{+})$ to (\ref{adjoint_DO1}) and obtain
\begin{eqnarray}
\sum_{m=1}^{N_T} \int_{I_m} \big(-(\partial_t p_{\delta}, \psi_{\delta})
&+& a_{h}^{s}(\psi_{\delta},p_{\delta})\big)  dt - \sum \limits_{m=1}^{N_T}([p_{\delta}]_{m}, \psi_{m,\delta}^- ) + (q_{\delta, N_{T}}^{-}, p_{\delta, N_T}^{+})\notag\\
&=& - \sum_{m=1}^{N_T} \int_{I_m}(y_{\delta}-y_{\delta}^{d},\psi_{\delta})\; dt, \quad \forall \psi_{\delta} \in V_{h}^{k, q}.\label{adjoint_DO2}
\end{eqnarray}
On each subinterval $I_m$, the adjoint equation reads as
$$
\int_{I_m} \left(-(\partial_t p_{\delta}, \psi_{\delta}) + a_{h}^{s}(\psi_{\delta},p_{\delta})\right)  dt - ([p_{\delta}]_{m}, \psi_{m,\delta}^- )
= -\int_{I_m}(y_{\delta}-y_{\delta}^{d},\psi_{\delta})\; dt.
$$
However, $(q_{\delta, N_{T}}^{-}, p_{\delta, N_T}^{+})$ does not  match the right-hand side of (\ref{adjoint_DO2}), so it is set to zero, i.e. $p_{\delta, N}^{+} = 0$. Now, we use  $a_{h}^{s}(\psi_{\delta},p_{\delta}) = a_{h}^{a}(p_{\delta},\psi_{\delta})$. Thus, we arrive at (\ref{fully_discerete_ocp_b}). Therefore, \emph{OD} and \emph{DO} approaches commute.
\section{Numerical Results}\label{numerical}
In this section, we present numerical results.  The state, the adjoint, and the control variables are discretized using the piecewise linear polynomials in space. The discretized control problem are solved by the primal dual active set (PDAS) algorithm \cite{MBergounioux_MHaddou_MHintermueller_KKunisch_2001}.
In order to measure the error in the state and adjoint approximation in terms of $L^{\infty}(0,1;L^{2}(\Omega))$ norm, the error in the control approximation in terms of $L^{2}(0,1;L^{2}(\Omega))$ norm. In all numerical examples, we have taken $h=k$.

We note that, in the case of dG(0) method, the approximating polynomials are piecewise constant in time and the resulting scheme is a version of the backward Euler method with a modified right-hand side \cite[Chapter 7]{VThomee_1997}
\begin{eqnarray*}
  (M + k A^s)y_{h,m} &=& M y_{h,m-1} + \frac{k}{2}(f_{h,m} + f_{h,m-1}) + \frac{k}{2}M(u_{h,m} + u_{h,m-1}),\\
  (M + k A^a)p_{h,m-1}&=& M p_{h,m} - \frac{k}{2}M(y_{h,m} + y_{h,m-1}) + \frac{k}{2}(y^{d}_{h,m} + y^{d}_{h,m-1}).
\end{eqnarray*}
For dG(1) method, we use piecewise linear polynomials in time. The resulting linear system for the state on each time step is given as follows \cite[Chapter 7]{VThomee_1997}:
\begin{equation} \label{state_dG1}
\left(
  \begin{array}{cc}
    M + k A^s & M + \frac{k}{2}A^s \\
    \frac{k}{2}A^s &  \frac{1}{2}M + \frac{k}{3} A^s \\
  \end{array}
\right)
\left(
  \begin{array}{c}
    Y_0 \\
    Y_1 \\
  \end{array}
\right) =
\left(
  \begin{array}{c}
    M y_{h,m-1} + \frac{k}{2}(f_{h,m} + f_{h,m-1}) + \frac{k}{2}M(u_{h,m} + u_{h,m-1}) \\
    \frac{k}{2}(f_{h,m} + Mu_{h,m} ) \\
  \end{array}
\right),
\end{equation}
where $A^s$, $M$ are the stiffness matrix of the state equation and the mass matrix, respectively. We derive the solution at the time step $t_m$ as $y_{h,m} = Y_0 + Y_1$.  For the adjoint equation, we have
\begin{equation} \label{adjoint_dG1}
\left(
  \begin{array}{cc}
    M + k A^a & M + \frac{k}{2}A^a \\
    \frac{k}{2}A^a &  \frac{1}{2}M + \frac{k}{3} A^a \\
  \end{array}
\right)
\left(
  \begin{array}{c}
    P_0 \\
    P_1 \\
  \end{array}
\right) =
\left(
  \begin{array}{c}
    M p_{h,m} - \frac{k}{2}M(y_{h,m} + y_{h,m-1}) + \frac{k}{2}(y^{d}_{h,m} + y^{d}_{h,m-1}) \\
    -\frac{k}{2}(M y_{h,m-1} - y^{d}_{h,m-1} )\\
  \end{array}
\right),
\end{equation}
where $A^{a}$ is the stiffness matrix for the adjoint equation. We obtain the adjoint at the time step $t_{m-1}$ as $p_{h,m-1} = P_0 + P_1$. We apply block-partitioning to these linear systems once in order to solve the systems for each time interval.\\

The main drawback of the dG time discretization is the solution of large coupled linear systems in block form. Several solvers are suggested to overcome this especially for nonlinear problems \cite{vexler13}. Because we are using constant time steps, the coupled matrices on the righthand side  of (\ref{state_dG1}) and (\ref{adjoint_dG1}) have to decomposed (LU block factorization) at the begin of the integration. Then the the state and adjoint equations are solved at each time step by forward elimination and back substitution using the block factorized matrices.  \\

\textbf{Example 1}: We consider the problem in \cite{HFu_HRui_2009a} with the following parameters by adding the reaction term
 \[
Q=(0,1] \times \Omega, \; \Omega=(0,1)^{2}, \; \epsilon=10^{-5}, \; \beta=(1,0)^T, \; r=1, \; \alpha=1 \text{ and} \;  u \ge 0.
\]
The source function $f$, the desired state $y_{d}$ and the initial condition $y_0$ are computed from (\ref{s4}) using the following exact solutions of the  state, adjoint and control, respectively,
\begin{eqnarray*}
y(x,t) &=& \exp(-t) \sin(2 \pi x_1) \sin(2 \pi x_2), \\
p(x,t) &=& \exp(-t) (1-t) \sin(2 \pi x_1) \sin(2 \pi x_2), \\
u(x,t) &=& \max\left(0,-\frac{1}{\alpha}p\right).
\end{eqnarray*}

In Table \ref{T:1}, errors and converge rates for dG(0) and backward Euler method are shown. For dG(0) and backward Euler method leads the first order convergence, due to the dominance of temporal errors, which is optimal in time.

\begin{table}
\centering
\caption{Example 1 by dG(0) and backward Euler(in parenthesis) method.}
\label{T:1}
\centering
\begin{tabular}{l c c c c c c }
\hline
$k$  & $\|y-y_{\delta}\|$ & Rate & $\|p-p_{\delta}\|$ & Rate & $\|u-u_{\delta}\|$ & Rate  \\
\hline
$\frac{1}{5}$   & 4.41e-2(2.45e-2) &-(-)      &8.77e-2(3.39e-2)&-(-)      &4.37e-2(2.43e-2)&-(-)\\
$\frac{1}{10}$  & 2.22e-2(6.84e-3) &0.99(1.84)&4.53e-2(1.42e-2)&0.95(1.26)&1.77e-2(5.73e-3)&1.31(2.08)\\
$\frac{1}{20}$  & 1.18e-2(6.84e-3) &0.99(1.53)&2.35e-2(6.88e-3)&0.95(1.04)&8.63e-3(2.67e-3)&1.03(1.10)\\
$\frac{1}{40}$  & 6.20e-3(6.84e-3) &0.93(1.17)&1.20e-2(3.45e-3)&0.96(1.00)&4.28e-3(1.34e-3)&1.01(1.00)\\
\hline
\end{tabular}
\end{table}

In Table \ref{T:2}, errors and converge rates for Crank-Nicolson method obtained by \emph{OD} and \emph{DO} approaches are shown. For Crank-Nicolson method, \emph{OD} approach  optimal second order convergence is achieved. But in the \emph{DO} approach, the discretization of the right-hand side of the adjoint by a one-step method is reflected affects the numerical results and the optimal order of convergence is not achieved.
\begin{table}[ht]
\centering
\caption{Example 1 by Crank-Nicolson method \emph{OD} and \emph{DO} approach(in parenthesis).}
\label{T:2}
\centering
\begin{tabular}{l c c c c c c }
\hline
$k$  & $\|y-y_{\delta}\|$ & Rate & $\|p-p_{\delta}\|$ & Rate & $\|u-u_{\delta}\|$ & Rate  \\
\hline
$\frac{1}{5}$   & 5.38e-2(5.31e-2) &-(-)      &3.22e-2(4.16e-1)&-(-)      &2.18e-2(4.33e-2)&-(-)\\
$\frac{1}{10}$  & 1.35e-2(1.36e-2) &1.99(1.97)&8.19e-3(1.90e-1)&1.98(1.13)&3.68e-3(1.24e-2)&2.57(1.80)\\
$\frac{1}{20}$  & 3.41e-3(3.43e-3) &1.99(1.98)&2.07e-3(9.10e-2)&1.98(1.06)&9.34e-4(4.38e-3)&1.98(1.50)\\
$\frac{1}{40}$  & 8.58e-4(8.65e-4) &1.99(1.99)&5.02e-4(4.45e-2)&2.05(1.03)&2.13e-4(1.63e-3)&2.13(1.42)\\
\hline
\end{tabular}
\end{table}

\begin{table}[ht]
\centering
\caption{Example 1 by dG(1) method.}
\label{T:3}
\centering
\begin{tabular}{l c c c c c c }
\hline
$k$  & $\|y-y_{\delta}\|$ & Rate & $\|p-p_{\delta}\|$ & Rate & $\|u-u_{\delta}\|$ & Rate  \\
\hline
$\frac{1}{5}$   & 3.65e-2 &-&5.36e-2&-&4.34e-2&-\\
$\frac{1}{10}$  & 8.59e-3 &2.09&1.35e-2&1.99&6.71e-3&2.70\\
$\frac{1}{20}$  & 2.14e-3 &2.00&3.35e-3&2.02&1.56e-3&2.10\\
$\frac{1}{40}$  & 5.36e-4 &2.00&8.16e-4&2.04&3.61e-4&2.11\\
\hline
\end{tabular}
\end{table}

In Table \ref{T:3}, results for dG(1) time discretization is shown and indicates that the second order convergence is achieved. The error in the state is smaller than for Crank-Nicolson method with \emph{OD} approach, while the errors in adjoint and the control are close for both discretizations. \\

\textbf{Example 2}: We consider the problem in \cite{HFu_2010a} with the following parameters by adding the reaction term
 \[
Q=(0,1] \times \Omega, \; \Omega=(0,1)^{2}, \; \epsilon=10^{-5}, \; \beta=(0.5, 0.5)^T, \; r=1, \; \alpha=1 \; \text{ and} \; 0 \leq u \leq 0.5.
\]
The source function $f$, the desired state $y_{d}$ and the initial condition $y_0$ are computed from (\ref{s4}) using the following exact solutions of the state, adjoint and control, respectively,

\begin{eqnarray*}
p(x,t) &=& \sin(\pi t) \sin(2 \pi x_1) \sin(2 \pi x_2) \exp\left( \frac{-1+\cos(t_x)}{\sqrt{\varepsilon}} \right) , \\
y(x,t) &=& p\left( \frac{1}{2 \sqrt{\varepsilon}}\sin(t_x) + 8 \varepsilon \pi^{2} + \frac{\sqrt{\varepsilon}}{2}\cos(t_x) - \frac{1}{2}\sin^{2}(t_x) \right) \\&-& \pi \cos(\pi t) \sin(2 \pi x_1)\sin(2 \pi x_2) \exp\left( \frac{-1+\cos(t_x)}{\sqrt{\varepsilon}}\right)\\
u(x,t) &=& \max \left(0,\min(-\frac{1}{\alpha}p, 0.5 )\right).
\end{eqnarray*}

\begin{figure}[htb]
  \centering
  \includegraphics[height=0.50\textwidth]{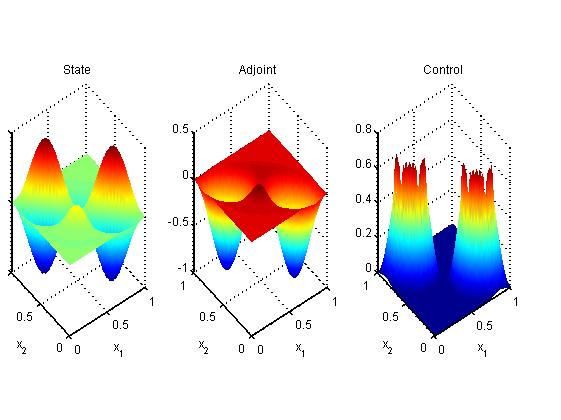}
   \caption{Example 2: Approximate solutions at t= 0.5 with Crank-Nicolson \emph{DO} approach.}
  \label{Fig:Ex2_CN_DO}
\end{figure}

\begin{figure}[htb]
  \centering \includegraphics[height=0.50\textwidth]{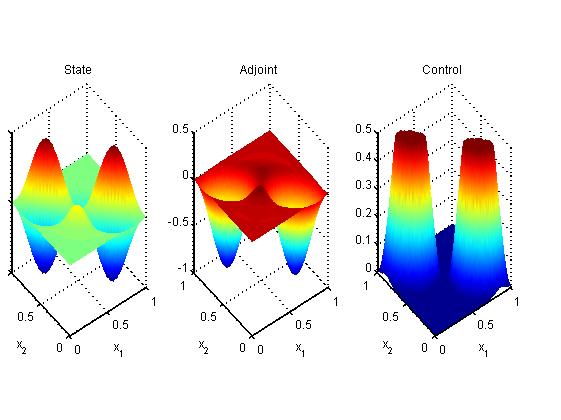}
   \caption{Example 2: Approximate solutions at t= 0.5 with Crank-Nicolson \emph{OD} approach.}
  \label{Fig:Ex2_CN_OD}
\end{figure}

As opposed to the previous example, the exact solution of PDE constrained depends on the diffusion explicitly and the problem is highly convection dominated. This example cannot be solved properly by using dG(0) and backward Euler method for time discretization. Therefore, we present errors for Crank-Nicolson method in Table \ref{T:4}, where the differences between \emph{OD} and \emph{DO} can be seen clearly. \emph{DO} approach causes order reduction for adjoint and control. We observe that \emph{DO} approach produces oscillations in the control Figure~\ref{Fig:Ex2_CN_DO}, whereas the approximate solutions for the control are smooth  for  \emph{OD} Figure~\ref{Fig:Ex2_CN_OD}. However, due to the convection dominated nature of the problem, the optimal order of convergence cannot be achieved in the \emph{OD} approach in contrast to the {\em Example 1}.

\begin{table}[ht]
\centering
\caption{Example 2 by Crank-Nicolson method \emph{OD} and \emph{DO} approach(in parenthesis).}
\label{T:4}
\centering
\begin{tabular}{l c c c c c c }
\hline
$k$  & $\|y-y_{\delta}\|$ & Rate & $\|p-p_{\delta}\|$ & Rate & $\|u-u_{\delta}\|$ & Rate  \\
\hline
$\frac{1}{5}$   & 2.32(2.32)  &    -(-)   &   3.17e-1(3.21e-1) &    -(-)      &    1.51e-1(1.43e-1) &   -(-)     \\
$\frac{1}{10}$  & 1.05(1.05)  & 1.14(1.14)&   1.25e-1(1.26e-1) &   1.34(1.35) &    5.09e-2(4.93e-2) &  1.57(1.54) \\
$\frac{1}{20}$  & 3.72e-1(3.74e-1) &  1.50(1.50)&   6.35e-2(7.47e-2) &   0.97(0.76) &    2.45e-2(3.36e-2) &  1.05(0.55) \\
$\frac{1}{40}$  & 1.09e-1(1.10e-1) &  1.77(1.76)&   2.17e-2(3.57e-2) &   1.55(1.07) &    8.31e-3(2.07e-2) &  1.56(0.70) \\
\hline
\end{tabular}
\end{table}


\begin{figure}[htb]
  \centering \includegraphics[height=0.50\textwidth]{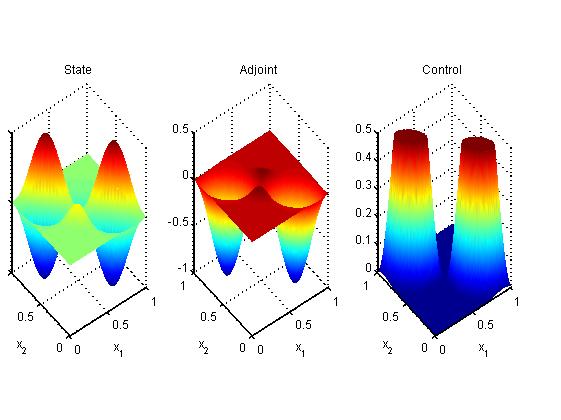}
\caption{Example 2: Exact solutions at t=0.5.}
   \label{Fig:Ex2_exact}
\end{figure}

In Table \ref{T:5}, numerical results for dG(1) discretization is shown. As opposed to the results in Table \ref{T:4}, the error in state, adjoint and control are smaller than in case of CN and the optimal quadratic convergence is achieved. 
\begin{table}[ht]
\caption{Example 2 by dG(1) method.}
\label{T:5}
\centering
\begin{tabular}{l c c c c c c }
\hline
$k$  & $\|y-y_{\delta}\|$ & Rate & $\|p-p_{\delta}\|$ & Rate & $\|u-u_{\delta}\|$ & Rate  \\
\hline
$\frac{1}{5}$   & 2.25e+0 &-&3.30e-1&-&1.48e-1&-\\
$\frac{1}{10}$  & 6.15e-1 &1.87&5.50e-2&2.58&2.38e-2&2.63\\
$\frac{1}{20}$  & 1.34e-1 &2.20&1.45e-2&1.92&8.01e-3&1.57\\
$\frac{1}{40}$  & 2.65e-2 &2.34&3.13e-3&2.22&2.27e-3&1.82\\
\hline
\end{tabular}
\end{table}

\begin{figure}[htb]
\centering  \includegraphics[height=0.50\textwidth]{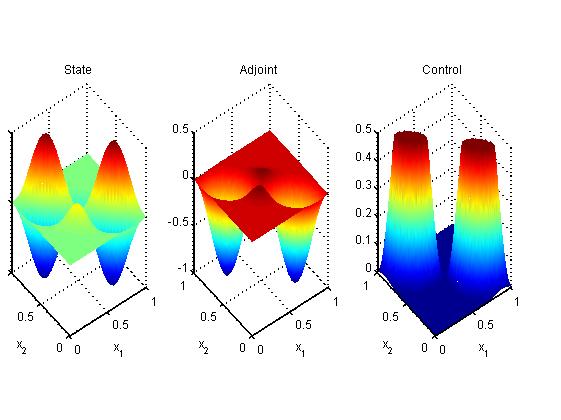}
\caption{Example 2: Approximate solutions at t=0.5 with dG(1) method.}
   \label{Fig:Ex2_dG1}
\end{figure}

In Figure~\ref{Fig:Ex2_dG1}, we present the exact and the approximate solution at $t = 0.5$ showing that the problem is approximated accurately.
\section{Conclusions}\label{conclusion}
For dG in time discretization, the numerical results confirm convergence rates and  \emph{DO}, \emph{OD} approaches commute. In a future work, we will study derivation of the optimal convergence rates under lower regularity assumptions and we will apply space-time adaptivity for convection dominated problems with boundary or interior layers.

\section*{Acknowledgement}
The authors thank to Konstantinos Chrysafinos for his explanations regarding error estimates and references. This research was supported by the Middle East Technical University Research Fund Project (BAP-07-05-2012-102).

\bibliographystyle{plain}

\begin{thebibliography}{00}%
%
  \bibitem{Apel_Flaig_2012_a}
  T. Apel, T.G. Flaig, {C}rank-{N}icolson schemes for optimal control problems with evolution equations, SIAM J. Numer. Anal. 50(3) (2012) 1482-1512.

\bibitem{DNArnold_FBrezzi_BCockburn_LDMarini_2001a}
D.N. Arnold, F. Brezzi, B. Cockburn, L.D. Marini, Unified analysis of discontinuous {G}alerkin methods for
              elliptic problems, SIAM J. Numer. Anal. 39(5) (2001/02) 1749-1779.

\bibitem{RBecker_BVexler_2007a}
 R. Becker, B. Vexler, Optimal control of the convection-diffusion equation using stabilized finite element
 methods, Numer. Math. 106(3) (2007) 349-367.

\bibitem{MBergounioux_MHaddou_MHintermueller_KKunisch_2001}
M. Bergounioux, M. Haddou, M. Hintermueller, K. Kunisch, A comparison of interior–-point methods and a Moreau–Yosida based active set strategy for constrained optimal control problems, SIAM J. Optim. 11(2) (2000) 495-521.

 \bibitem{EBurman_2011a}
   E. Burman, {C}rank-{N}icolson finite element methods using symmetric stabilization with an application to optimal control problems subject to transient advection-diffusion equations, Comm. Math Sci. 9(1) (2011) 319-329.


 \bibitem{CKonstantinos_2007a}
   K. Chrysafinos, Discontinuous {G}alerkin approximations for distributed optimal control problems constrained by parabolic {PDE}'s,
   Int. J. Numer. Anal. Model. 4(3-4) (2007) 690-712.


  \bibitem{KChrysafinos_NJWalkington_2006b}
   K. Chrysafinos, N.J. Walkington, Error estimates for the discontinuous {G}alerkin methods for parabolic equations, SIAM J. Numer. Anal. 44(1) (2006) 349-366.

\bibitem{SSCollis_MHeinkenschloss_2002a}
S.S. Collis, M. Heinkenschloss, Analysis of the streamline upwind/Petrov Galerkin method applied
to the solution of optimal control problems. Tech. Rep. TR02–01, Department of Computational and
Applied Mathematics, Rice University, Houston, TX 77005-1892 (2002).



   \bibitem{KEriksson_CJohnson_VThomee_1985}
  K. Eriksson, C. Johnson, V. Thom{\'e}e, Time discretization of parabolic problems by the discontinuous
  {G}alerkin method, RAIRO Mod\'el. Math. Anal. Num\'er. 19(4) (1985) 611-643.

 \bibitem{MFeistauer_VKucera_KNajzar_JProkopovza_2011a}
 M. Feistauer, V. Ku{\v{c}}era, K. Najzar, J. Prokopov{\'a}, Analysis of space-time discontinuos Galerkin method for nonlinear convection-diffusion problems, Numer. Math. 117(2) (2011) 251-288.


\bibitem{HFu_2010a}
   H. Fu, A characteristic finite element method for optimal control problems
governed by convection-diffusion equations, J. Comput. Appl. Math. 235 (2010) 825-836.

   \bibitem{HFu_HRui_2009a}
   H. Fu, H. Rui, A priori error estimates for optimal control problems governed by transient advection-diffusion
   equations, J. Sci. Comput. 38(3) (2009) 290-315.


 \bibitem{MHinze_NYan_ZZhou_2009a}
 M. Hinze, N. Yan, Z. Zhou, Variational discretization for optimal control governed by convection
 dominated diffusion equations, J. Comp. Math. 27(2-3) (2009) 237-253.

 \bibitem{DLeykekhman_2012b}
 D. Leykekhman, Investigation of commutative properties of discontinuous Galerkin methods in PDE
 constrained optimal control problems, J. Sci. Comput. 53(3) (2012) 483-511.

 \bibitem{DLeykekhman_MHeinkenschloss_2012a}
 D. Leykekhman, M. Heinkenschloss, Local error analysis of discontinuous Galerkin methods for
 advection-dominated elliptic linear-quadratic optimal control problems, SIAM J. Numer. Anal. 50(4) (2012) 2012-2038.

\bibitem{JLLions_1971}
J.L. Lions, Optimal Control of Systems Governed by Partial Differential Equations, Springer Verlag,
Berlin, Heidelberg, New York, 1971.

 \bibitem{DMeidner_BVexler_2008b}
 D. Meidner, B. Vexler, A priori error estimates for space-time finite element discretization of parabolic
 optimal control problems. II. Problems with control constraints, SIAM J. Control Optim. 47(3) (2008), 1301-1329.

\bibitem{BRiviere_2008a}
B. Riv\`{\i}ere, Discontinuous Galerkin Methods for Solving Elliptic and Parabolic Equations: Theory and
Implementation, Frontiers in Applied Mathematics, vol. 35. SIAM, Philadelphia, 2008.


\bibitem{DSchotzau_LZhu_2009a}
D. Sch{\"o}tzau, L. Zhu, A robust a-posteriori error estimator for discontinuous {G}alerkin methods for
convection-diffusion equations, Appl. Numer. Math. 59(9) (2009) 2236-2255.


\bibitem{MStoll_AWathen_2010a}
M., Stoll, A., Wathen, A.: All-at-once solution of time-dependent PDE-constrained optimization problems. Tech. Rep. TR2, Max Planck Institute for Dynamics of Complex Technical Systems, 39106, Magdeburg (2010)

   \bibitem{TSun_2010a}
   T, Sun, Discontinuous Galerkin finite element method with interior penalties for convection diffusion
   optimal control problem, Int. J. Numer. Anal. Model. 7(1) (2010) 87-107.

 \bibitem{VThomee_1997}
 V. Thom\'{e}e, {G}alerkin Finite Element Methods for Parabolic Problems, 2nd Ed., Springer Verlag, Berlin, 2006.

\bibitem{FTroltzsch_2010a}
F. Tr\"{o}ltzsch, Optimal Control of Partial Differential Equations: Theory,
Methods and Applications, Graduate Studies in Mathematics, American
Mathematical Society, 112, Providence, RI, 2010.


\bibitem{vexler13} Richter, T., Springer, A., and Vexler, B.,
Efficient numerical realization of discontinuous
Galerkin methods for temporal discretization. 124 (2013) 151–182.

  \bibitem{MVlasak_VDolejsi_JHajek_2010a}
   M. Vlas{\'a}k, V. Dolej{\v{s}}{\'{\i}}, J. H{\'a}jek, A priori error estimates of an extrapolated space-time
   discontinuous {G}alerkin method for nonlinear convection-diffusion problems, Numer. Methods Partial Differential Equations. 27(6) (2011) 1456-1482.

\bibitem{HYucel_MHeinkenschloss_BKarasozen_2012b}
H. Y\"{u}cel, M. Heinkenschloss, B. Karas\"{o}zen, Distributed optimal control of diffusion-convection-reaction equations using discontinuous Galerkin methods, Proceedings of ENUMATH 2011, Springer, Berlin, 389-397 (2013).

\bibitem{HYucel_BKarasozen_2013}
Y\"{u}cel, H., Karas\"{o}zen, B.: Adaptive Symmetric Interior Penalty Galerkin (SIPG) method for
optimal control of convection diffusion equations with control constraints. Optimization, (2013).

   \bibitem{ZZhou_NYan_2010a}
   Z. Zhou, N. Yan, The local discontinuous Galerkin method for optimal control problem governed by
   convection diffusion equations, Int. J. Numer. Anal. Model. 7(4) (2010) 681-699.

\end{thebibliography}

\end{document}